\documentclass[12pt]{article}%
\usepackage[applemac]{inputenc}
\usepackage{graphicx}  
\usepackage{caption}   
\usepackage{color, colortbl}
\usepackage{booktabs}
\usepackage{amsmath,amssymb}
\usepackage[T1]{fontenc}
\usepackage[T1]{fontenc}
\usepackage{amsfonts}
\usepackage[applemac]{inputenc}
\usepackage{amsmath,amssymb}
\usepackage{color}
\usepackage{color, colortbl}
\usepackage{amsmath}
\usepackage{amssymb}
\usepackage{graphicx}
\usepackage{tikz}
\usepackage{geometry}
\usetikzlibrary{arrows}
\usepackage{amsfonts}
\usepackage{amsmath,amssymb}
\usepackage[applemac]{inputenc}
\usepackage{color}
\usepackage{amsmath}
\usepackage{amssymb}
\usepackage{graphicx}
\usepackage{tikz}
\newtheorem{theorem}{Theorem}[section]

\newtheorem{lemma}[theorem]{Lemma}
\newcommand{\E}{E}
\newtheorem{proposition}[theorem]{Proposition}
\newtheorem{corollary}[theorem]{Corollary}
\newtheorem{definition}[theorem]{Definition}

\newtheorem{example}[theorem]{Example}

\newtheorem{remark}[theorem]{Remark}

\usetikzlibrary{arrows,shapes,automata,petri}
\newenvironment{myenumerate}{

\begin{enumerate}}{\end{enumerate}}
\newcommand{\dproof}{\noindent {Proof.} \quad}
\newcommand{\fproof}{\hfill $\square$ \bigskip}

\numberwithin{equation}{section}

\definecolor{LightCyan}{rgb}{0.88,1,1}
\def\cF{{\mathcal {F}}}

\def\RR{{\mathbb{ R}}}
\def\cB{{\mathcal{B}}}

\def\1B{\text{1\!\!I}}

\def\l{\langle}
\def\<{\langle}
\def\>{\rangle}
\def\P{\mathbb{P}}
\def\R{\mathbb{R}}
\def\l{\lambda}
\def\S{\mathcal{S}}
\def\E{\mathbb{E}}
\def\N{\mathbb{N}}

\def\F{\mathcal{F}}

\def\realio{\mathbb{R}_{0}}

\def\d{\diamond}

\def\d{\diamond}

\begin{document}

\title{Multiparameter L\'evy white noise theory and applications}

\author{Olfa Draouil$^{1}$, Rahma Yasmina Moulay Hachemi$^{2}$ \& Bernt \O ksendal$^{3}$
}
\date{26 September 2026}
\maketitle 
\footnotetext[1]{Department of Mathematics, University of Tunis El Manar, Tunisia\newline
Email: olfa.draouil@fst.utm.tn}

\footnotetext[2]{%
Department of Mathematics, University of Oslo, Norway.\newline
Email: ryhachem@math.uio.no}

\footnotetext[3]{%
Department of Mathematics, University of Oslo, Norway. \\
Email: oksendal@math.uio.no.}


\begin{abstract}
We construct a white noise theory and white noise calculus for the (multi-parameter) L\' evy sheet and its compensated Poisson random measures. The theory is applied to a stochastic partial differential equation subject to L\' evy noise.
  \end{abstract}

\section{Introduction}

The theory of white noise for Brownian motion was initiated by T.Hida \cite{H} and further developed by him and coauthors Kuo, Potthoff and Streit in \cite{HKPS}. Subsequently an extension to (multiparameter)  Brownian sheet was carried out and applied to stochastic partial differential equations by Holden, Ub\o e, \O ksendal and Zhang  in \cite {HOUZ}. 

Later Di Nunno, \O ksendal \& Proske \cite{DOP1} and \O ksendal \& Proske \cite{OP} developed a white noise theory for L\' evy processes and for compensated Poisson random measures, respectively. See the book \cite{DOP} for a comprehensive coverage of this.  

The theory of (multiparameter) L\' evy sheets was developed by Iafrate \& Ricciutu  in \cite{IR} and by Dalang \& Walsh \cite{DW}. See also the earlier paper Adler et al  \cite{AMSW}. 

Subsequently a white noise theory for L\' evy sheets was carried out in Chapter 5 of the book Holden et al \cite{HOUZ}. 

The purpose of this paper is to develop this theory further, in particular we prove a fundamental relation between Wick products and
Skorohod integrals. Moreover, we apply this theory to study a stochastic  equation modeling population growth in a random environment, driven by both Brownian and
L\' evy multiparameter white noise.

The content of our paper is the following:
\begin{itemize}
    \item 
    In Section 2 we recall the theory of L\' evy sheets and we give a construction of L\' evy sheets based on the Bochner-Minlos theorem in the white noise probability space.
    \item 
    In Section 3 we present a white noise theory for multiparameter L\' evy sheets.
    \item 
    In Section 4 we develop a white noise calculus based on the Wick product
    \item 
    Finally, in Section 5 we apply the theory to study a time-fractional stochastic partial differential equation modeling population growth in a Brownian-L\' evy noisy environment.
\end{itemize}

\section{Multiparameter L\' evy sheets and L\' evy White Noise}


\subsection{Multiparameter L\'{e}vy sheets}

We recall the theory of L\' evy sheets, mainly based on \cite{IR}, \cite{DW}:

\begin{definition}
  A \emph{multiparameter L\'{e}vy sheet} is a stochastic process $L(x)=L(x,\omega); x \in \R_+^d,\omega \in \Omega$ defined on a probability space $(\Omega,\F,\P)$ and satisfying the following conditions:
\begin{itemize}
\item
    $L$(0)=0 a.s.
\item
$L$ has stationary, independent increments with respect to the natural partial ordering on $\R_+^d$
\item
$L$ is continuous in probability and has c\'adl\'ag paths a.s. 
\end{itemize}
\end{definition}

\begin{remark}
If $R=\Pi_{j=1}^n (x_j,x'_j]$ is a box in $\R_+^d$ then $\Delta_R L$ denotes the increment of $L$ over the box $R$. For example, if $d=2$ then
$$\Delta_R L= L(x_1',x_2')-L(x_1',x_2) - L(x_1, x_2') + L(x_1,x_2).$$
\end{remark}
\begin{remark}    
If we also assume that $L$ has continuous paths, then $L$ becomes a \emph{Brownian sheet}.
\end{remark}

\begin{definition}
\begin{myenumerate}
    \item
    The jump of $L$ at the point $x \in \R_{+}^d$ is defined by
    \begin{align}
        \Delta L_x= L(x)-L(x^{-}).
    \end{align}
    \item
Let $R$ denote a closed finite box in $\R_+^d$ not containing 0, and let $U$ be a Borel subset of $\R_0:=\R \setminus 0$. Note that because the paths are c\'adl\'ag there are only finitely many jumps in $R$
Define \emph{the jump measure} $N(R,U)$ to be the number of jumps of $L$ within the box $R$ of jumps size in $U$. By stationarity the map $R \mapsto N(R,U);  U \in \mathcal{B}(\R_0)$ (the family of all Borel subsets of $\R_0$) depends only on the volume of $R$ and we can write the random measure generated by $(R,U) \mapsto N(R,U)$ in differential form as  $N(dx, du)$, where $dx$ and $du$ are the Lebesgue measure differentials corresponding to $R$ and $U$, respectively.
\item
Put $R_1=[0,1]^d$. The measure $U \mapsto \nu(U):=\E [N(R_1,U)];  U \in \mathcal{B}(\R_0)$ is called the L\'{e}vy measure of $L$.
\item
The corresponding \emph{compensated jump measure} is defined by
$$ \widetilde{N}(dx,dz)= N(dx,dz) - \nu(dz) dx,$$
where $dx=dx_1 dx_2 ... dx_n$ is the volume differential at $x=(x_1, x_2, ... ,x_d) \in \R^d$.
\end{myenumerate}
\end{definition}
\subsection{The L\'{e}vy white noise probability space and the pure jump L\'{e}vy sheet}

In the following we let $\S=\S(\R^d)$ denote the Schwartz space of rapidly decreasing smooth functions on $\R^d$, and we let $\S'=\S'(\R^d))$ denote its dual, which is the space of tempered distributions on $\R^d$, equipped with the weak$^{*}$ topology. 

It is possible to construct a \emph{pure jump L\'{e}vy sheet} $L$ corresponding to a given  L\'{e}vy measure $\nu$  
as a process on the Borel $\sigma$-algebra $\cF = \cB(\S'(\R^d))$. The details are the following:

Let $\nu$ be a measure on $\R_0$ such that
\begin{align}
\int_{\R_0} z^2 \nu(dz) < \infty.
\end{align}
By the Bochner-Minlos-Sazonov theorem\index{Bochner-Minlos-Sazonov theorem} (see \cite{GV}), there exists a probability measure $\P$ on the Borel $\sigma$-algebra $\cF = \cB(\S'(\R^d))$ such that
\begin{equation}
\int_\Omega e^{i\langle\omega,\varphi \rangle}\P(d\omega)=\exp\left(\int_{\R^n}\Psi(\varphi(x))dx\right);
\quad \varphi\in\mathcal{S}(\R^d), \label{boch}
\end{equation}
where $i=\sqrt{-1}$ and $\langle\omega,\varphi\rangle$ denotes the action of $\omega\in \Omega:=\mathcal{S}'(\R^d)$ on $\varphi\in \mathcal{S}(\R^d)$ and 
\begin{equation}
\Psi(w)=\int_{\R}\left(e^{i\,w\,z}-1-i\,w\,z\right)\nu(dz); w \in \R.
\end{equation}

We let $\E[\cdot]=\E_{\P}[\cdot]$ and $Var[\cdot]=Var_{\P}[\cdot]$ denote expectation and variance, respectively, with respect to $\P$.\\
Note that by \eqref{boch} we have, for all parameters $a \in \R$,
\begin{align}
    \E| e^{i\langle \omega,a \varphi \rangle)}]=\exp\left(\int_{\R^n}\int_{\R_0} \{ e^{i a \varphi(x)z} - 1 - i a \varphi(x)z\} \nu(dz) dx \right) \label{boch2}
\end{align}
Expanding the exponential function on both sides of \eqref{boch2} and comparing the terms with the same power of $a$, we obtain the following result:
\begin{lemma}
\label{2.5}
Let $\varphi\in\mathcal{S}(\R^d)$. 
Then we have
\begin{equation}
\E[\langle\cdot,\varphi\rangle]=0
\end{equation}
and
\begin{equation}
Var[\langle\cdot,\varphi \rangle]:=\E[\langle\cdot,\varphi\rangle^2] = M \int_{\R^d} \varphi^2(x)dx, \label{Var}
\end{equation}
where
\begin{align}
    M=\int_{\R_0} z^2 \nu(dz).
\end{align}
Therefore,
\begin{align}
    \E[(\langle \omega,\varphi_1\rangle - \langle\omega,\varphi_2\rangle)^2]=M \int_{\R^d} (\varphi_1(x)-\varphi_2(x))^2 dx=M \| \varphi_1 - \varphi_2 \|^2_{L^2(\lambda)}, 
\end{align}
where $\lambda$ is the Lebesgue measure on $\R^d$.

\end{lemma}

Using this result, we can for each  $\omega \in \Omega$ extend the definition of $\langle \omega,f\rangle$ from $f\in\mathcal{S}(\R^d)$ to any $f\in L^2(\lambda)$ as follows:

If $f\in L^2(\lambda)$, choose $f_n\in\mathcal{S}(\R^n)$ such that $f_n\rightarrow f$ in $L^2(\lambda)$. 
Then by Lemma \ref{2.5} we see that $\{\langle
\omega,f_n\rangle\}_{n=1}^\infty$ is a Cauchy sequence in $L^2(\P)$ and hence convergent in $L^2(\P)$. 
Moreover, the limit depends only on
$f$ and not the sequence $\{f_n\}_{n=1}^\infty$. 
We denote this limit by $\langle \omega,f\rangle; \  \omega \in \Omega, f \in L^2(\lambda).$

In particular, if we define
$$L(x):=L(x,\omega):=\langle\omega,\chi_
{\Pi_{j=1}^d [0,x_j]} \rangle;\quad x=(x_1, x_2, ... ,x_d) \in \R_+^d,$$
then we have the following result:

\begin{theorem}
The stochastic process $L(x)$, $x \in \R_+^d$ has a c\`{a}dl\`{a}g version, also denoted by $L$. 
This process $L(x)\ ; x\in \R_+^d$ is a (pure jump) L\'{e}vy
sheet with L\'{e}vy measure $\nu$.
\end{theorem}

\section{Multiparameter L\' evy white noise theory}
Our presentation in this section is based on a combination of the results in the papers \cite{IR}, \cite{DW}, \cite{DOP} and \cite{OP} and the book \cite{HOUZ}.\\

In the following we will for simplicity assume that
\begin{align}
\int_{\R_0} z^2 \nu(dz) < \infty.
\end{align}
Then by the L\'{e}vy-Khinchine theorem the L\'{e}vy sheet $L$ can be given the following representation:
\begin{align}
    L(x)=\alpha \Pi_{j=1}^n x_j + \sigma B(x) + \int_{\R_0} z \widetilde{N}(x, dz); \ x \in \R_+^d,
\end{align}
for some constants $\alpha, \sigma$, where $B(x)$ is a Brownian sheet.

We will from now on assume that $\alpha=\sigma=0$. 
Then we have a \emph{pure jump L\'evy sheet}, with representation
\begin{align} 
L(x)=\int_{\R_0} z \widetilde{N}(x,dz)=\int_0^x \int_{\R_0} z \widetilde{N}(d\xi,dz); \ x \in \R_+^d. \label{LN}
\end{align}
In differential form this becomes
\begin{align}\label{ldx}
    L(dx)= \int_{\R_0} z \widetilde{N}(dx,dz).
\end{align}

We proceed to define the white noise of the L\'{e}vy sheet we constructed in the previous section. We will follow the presentation in \cite{DOP} Section 13.2  closely, but with the necessary modifications due to the fact that we are now dealing with an $n$-parameter L\'{e}vy sheet, not a classical 1-parameter L\'{e}vy process. Here are the details:\\

From now on we assume that the L\'{e}vy measure $\nu $ satisfies the
following condition:\\
For all $\varepsilon >0$ there exists $\lambda >0$ such that
\begin{equation}
\int_{\realio\backslash (-\varepsilon ,\varepsilon )}\exp (\lambda
\left\vert z\right\vert )\nu (dz)<\infty   \label{nu}
\end{equation}
One can prove that if this holds, then $\nu $ has finite moments of all orders $n\geq 2$. 
For example, the condition clearly holds if $\nu $ is supported on a compact subset of $\R$.

From this condition one can prove that the polynomials are dense in $L^{2}(\rho )$,
where
\begin{equation}
\rho (dz)=z^{2}\nu (dz).  
\end{equation}
See \cite{NualartSchoutens}. We let $\left\{\eta_{m}\right\} _{m\geq 0}=\left\{ \eta_0 =1,\eta_{1},\eta_{2},...\right\} $ \label{simb-1303} be the orthogonolization of
$\left\{1,z,z^{2},...\right\} $ with respect to the inner product of $L^{2}(\rho )$.

Define
\begin{equation}
p_{j}(z):=\left\Vert \eta_{j-1}\right\Vert _{L^{2}(\rho )}^{-1}z\eta_{j-1}(z);
\text{ }j=1,2,...  \label{10.3}
\end{equation}
and
\begin{equation}
m_{2}:=\left( \int_{\realio} z^{2}\nu (dz)\right) ^{\frac{1}{2}}
=\left\Vert \eta_{0}\right\Vert _{L^{2}(\rho )}
=\left\Vert 1\right\Vert_{L^{2}(\rho )}.  \label{10.4}
\end{equation}
In particular,
\begin{equation}
p_{1}(z)=m_{2}^{-1}z\text{ or }z=m_{2}p_{1}(z).  \label{10.5}
\end{equation}
Then $\left\{ p_{j}(z)\right\} _{j\geq 1}$ is an \emph{orthonormal basis}
for $L^{2}(\rho )$.

Recall that the {\it Hermite polynomials\/} $h_n(x); x \in \R$ are defined by
$$h_n(x)=(-1)^n
e^{{1}/{2}x^2}\frac{d^n}{dx^n}(e^{-{1}/{2}x^2});\quad
n=0,1,2,\cdots.\leqno{(2.2.1)}$$
Thus the first Hermite polynomials are
\begin{align}
h_0(x)&=1,\;h_1(x)=x,\;h_2(x)=x^2-1,\;h_3(x)=x^3-3x,\cr
h_4(x)&=x^4-6x^2+3,\;h_5(x)=x^5-10x^3+15x,\cdots.
\end{align}
The \emph{Hermite functions} $\xi_n(x)$ are defined by
$$\xi_n(x)=\pi^{-{1}/{4}}((n-1)!)^{-{1}/{2}}e^{-{1}/{2}x^2}h_{n-1}
(\sqrt{2}x);\quad n=1,2,\cdots.$$
We note the following  useful properties of $\xi_n$:
\begin{itemize}
\item
$\xi_n\in{\cal S}(\R)\quad\text{for all}\;n.$
\item
 The collection $\{\xi_n\}^\infty_{n=1}$ constitutes an
orthonormal basis for $L^2(\R).$
\item
$\sup\limits_{x\in\R}|\xi_n(x)|=\mathcal{O}(n^{-{1}/{12}}).$
\end{itemize}
See Hille and
Phillips (1957), Chapter 21, for proofs of these statements.
\vskip 0.2cm
Using these functions we define an orthogonal basis for $L^2(\P)$ as follows:\\
Let
$\beta=(\beta_1,\cdots,\beta_d)$ denote $d$-dimensional multi-indices with
$\beta_1,\cdots,\beta_d\in\N$, where $\mathbb{N}$ denotes the set of all natural numbers. By the above it follows that the family of tensor
products
$$\xi_\beta:=\xi_{(\beta_1,\cdots,\beta_d)}:=\xi_{\beta_1}\otimes\cdots
\otimes
\xi_{\beta_d}\;;\;\beta\in\N^d$$
forms an orthonormal basis for $L^2(\R^d)$. Let
$\beta^{(j)}=(\beta_1^{(j)},\beta_2^{(j)},\cdots,\beta_d^{(j)})$ be
the $j$th multi-index in some fixed ordering of all $d$-dimensional
multi-indices $\beta=(\beta_1,\cdots,\beta_d)\in\N^d$. We 
may assume that this ordering has the property that
$$i<j\Rightarrow\beta_1^{(i)}+\beta_2^{(i)}+\cdots+\beta_d^{(i)}\leq
\beta_1^{(j)}+\beta_2^{(j)}+\cdots+\beta_d^{(j)},$$
i.e., that the $\{\beta^{(j)}\}^\infty_{j=1}$ occur in increasing order.

Next, define the bijection $\kappa :\mathbb{N}\times \mathbb{N}\longrightarrow
\mathbb{N}$\label{simb-new2} by
\begin{equation}
\kappa (i,j)=j+(i+j-2)(i+j-1)/2;\quad  i,j=1,2, ...
\end{equation}
For $x=(x_1, x_2, ... , x_d) \in \R^d$ we proceed as in in Section 2.2.1 in \cite{HOUZ} and define $\left\{ e _{j}(x)\right\} _{j\geq 1}$ to be a fixed ordering of all $d$-times tensor products of the Hermite functions.
This family forms an orthonormal basis of $L^2(\R^d).$
Put
\begin{equation}
\theta _{\kappa (j,k)}(x,z):=e _{j}(x)p_{k}(z); \  x\in \R^n, z \in \R_0.
\end{equation}

In the following we let $\mathcal{J}$ denote the family of all finite multi-indices $\alpha=(\alpha_1, \alpha_2, ..., \alpha_m); \ m=1,2, ...$ of nonnegative integers $\alpha_j$. 

If $\alpha \in \mathcal{J}$ with $Index(\alpha )=j$  (the index of the last non-zero element) and $\left\vert \alpha
\right\vert =m$, we define $\theta ^{\otimes \alpha }: (\R^d \times \R_0)^m \mapsto \R$ by
\begin{eqnarray}
&&\theta ^{\otimes \alpha }(x^{(1)},z_{1};...;x^{(m)},z_{m}) \\
&=&\theta _{1}^{\otimes \alpha _{1}}\otimes ...\otimes \theta _{j}^{\otimes
\alpha _{j}}(x^{(1)},z_{1};... ;x^{(m)},z_{m})  \nonumber \\
&=&\underset{\alpha _{1}\text{ factors}}{\underbrace{\theta
_{1}(x^{(1)},z_{1})\cdot ...\cdot \theta _{1}(x^{(\alpha _{1})},z_{\alpha _{1}})}
}\cdot ...\cdot \underset{\alpha _{j}\text{ factors}}{\underbrace{\theta
_{j}(x^{(m-\alpha _{j}+1)},z_{m-\alpha _{j}+1})\cdot ...\cdot \theta
_{j}(x^{(m)},z_{m})}}.  \nonumber
\end{eqnarray}
We put $\theta _{i}^{\otimes 0}=1.$\\
Finally we let $\theta ^{\hat{\otimes }\alpha }$ denote the
\emph{symmetrized} tensor product\index{tensor product!symmetrized} of the $\theta _{k}$ $^{\prime }s:$
\begin{equation}
\theta ^{\hat{\otimes }\alpha }(x^{(1)},z_{1};...;x^{(m)},z_{m})
=\theta_{1}^{\hat{\otimes }\alpha _{1}}\otimes ...\otimes
\theta _{j}^{\hat{\otimes }\alpha _{j}}(x^{(1)},z_{1};...;x^{(m)},z_{m}).  
\end{equation}
The symmetrization is with respect to the $m$ pairs $(x^{(1)},z_1), (x^{(2)}z_2), \cdots (x^{(m)}z_m)$.
We set $K_0=1$ and for general $\alpha \in \mathcal{J}$ we define
\begin{equation}
K_{\alpha }:= I_{\left\vert \alpha \right\vert }\left(
\theta ^{\hat{\otimes }\alpha }\right),
\end{equation}
where in general $I_m$ is the symmetrized iterated integral with respect to $\widetilde{N}$. This integral is defined for symmetric functions $g: (\R^d \times \R_0)^m \mapsto \R$ by
\begin{align}
   I_m(g) &= \int_{(\R^d \times \R_0)^m } g(x^{(1)},z_1; x^{(2)},z_2; ... ;x^{(m)},z_m) \widetilde{N}^{\otimes m}(dx^{(1)} dz_1; ..., dx^{(m)} dz_m)\nonumber\\
   &:= (m!)^{d} J_m(g), 
\end{align}
where $J_m(g)$ is the $m$ times iterated integral
\begin{align}
    &J_m(g)\nonumber\\
    &=\int_{\R^d}\int_{\R_0} \Big(\int_{-\infty}^{x^{(m)}} \int_{\R_0}\Big( ... \Big(\int_{-\infty}^{x^{(2)}}g(x^{(1)},z_1; x^{(2)},z_2\cdots \nonumber\\
    &\cdots x^{(m)},z_m) ) \widetilde{N}(dx^{(1)},dz_1)\Big)  ... \Big) \widetilde{N}(dx^{(m)},dz_m). 
\end{align}
Here we are using the following notation, for $x=(x_1,x_2, ... ,x_m) \in \R^m$,
\begin{align}
    \int_{-\infty}^x= \int_{-\infty}^{x_m} \int_{-\infty}^{x_{m-1}}  ... \int_{-\infty}^{x_1}. 
\end{align}

\begin{example}
    For example, if $m=1$ we get
    \begin{align}
        I_1(g)=\int_{\R^d} \int_{\R_0} g(x,z) \widetilde{N}(dx,dz),
    \end{align}
    and if $m=2$ we have
    \begin{align}
        I_2(g)&=\int_{\R^d} \int_{\R_0} \Big( \int_{\R^d} \int_{\R_0} g(x^{(1)},z_1; x^{(2)},z_2) \widetilde{N}(dx^{(2)},dz_2) \Big)\widetilde{N}(dx^{(1)},dz_1)\nonumber\\
        &=2^d\int_{\R^d}\int_{\R_0} \Big(\int_{-\infty}^{x^{(2)}} \int_{\R_0} g(x^{(1)},z_1; x^{(2)},z_2) \widetilde{N}(dx^{(1)},dz_1)\Big) \widetilde{N}(dx^{(2)},dz_2)
    \end{align}
\end{example}
\begin{example}
\label{Ex10.1}
With $\varepsilon ^{(k)}=(0,...,0,1,0,...)$ with $1$ on $k$ th
place, we have, writing $\varepsilon ^{(\kappa (i,j))}=\varepsilon ^{(i,j)}$,
\begin{equation}
K_{\varepsilon ^{(i,j)}}
= I_1(\theta^{\otimes \varepsilon^{(i,j)}})
=I_{1}\left( \theta_{\kappa (i,j)} \right)
=I_{1}\left( e _{i}(x)p_{j}(z)\right) .  \label{10.11}
\end{equation}%

As in the Brownian sheet case one can prove that $\{ K_\alpha  \}_{\alpha \in \mathcal{J}}$ are
orthogonal in $L^2(\P)$ and
\[
\Vert K_\alpha \Vert^2_{L^2(\P)} = \alpha!.
\]
Note that if $|\alpha|=m$,
\[
\alpha ! =\Vert K_\alpha \Vert^2_{L^2(\P)} = m! \Vert \theta^{\hat\otimes \alpha} \Vert^2_{L^2((\lambda\times\nu)^m)},
\]
where as before $\lambda$ denotes Lebesgue measure on $\R^d$.
By our construction of $\theta ^{\hat{\otimes }\alpha }$ we know that
any $f\in \widetilde{L}^{2}((\lambda \times \nu )^{m})$ can be written%
\begin{equation}
f(x^{(1)},z_{1};...; x^{(m)},z_{m})=\sum_{\left\vert \alpha \right\vert
=m}c_{\alpha }\theta ^{\hat{\otimes }\alpha
}(x^{(1)},z_{1},..., x^{(m)},z_{m}).
\end{equation}
Hence
\begin{equation}
I_{m}(f_{m})=\sum_{\left\vert \alpha \right\vert =m}c_{\alpha }K_{\alpha }
\end{equation}
\end{example}
This gives the following chaos expansion:
\begin{theorem}
\label{Th10.2}
{\bf (Chaos expansion)}
\newline Any $F\in L^{2}(\P)$
has a unique expansion of the form
\begin{equation}
F=\sum_{\alpha \in \mathcal{J}}c_{\alpha }K_{\alpha }= \sum_{m=0}^\infty I_m(f_m).
\label{10.14}
\end{equation}
with
\begin{equation}
f_m=f_m (x^{(1)},z_{1};\cdots ; x^{(m)},z_{m})=\sum_{\left\vert \alpha \right\vert
=m}c_{\alpha }\theta ^{\hat{\otimes }\alpha
}(x^{(1)},z_{1},\cdots ; x^{(m)},z_{m}).
\end{equation}
and $c_{\alpha }\in \mathbb{R}$. 
Moreover,
\begin{equation}
\left\Vert F\right\Vert _{L^{2}(\P)}^{2}=\sum_{\alpha \in \mathcal{J}}\alpha
!c_{\alpha }^{2}.  \label{10.15}
\end{equation}
\end{theorem}

\begin{example}
\label{Ex10.3}
Let $h\in L^{2}(\mathbb{R})$ and define, for $x \in \R_+^d$,
\begin{align}
F=\int_0^{x} h(\xi)L(d\xi).
\end{align}
By  \eqref{LN} this can be written
\begin{align}
    F=I_{1}(h(\xi)\chi_{\Pi_{i=1}^d [0,x_i]}(\xi)z).
\end{align}
\end{example}
Therefore $F$ has the following expansion, writing $\chi_{\Pi_{i=1}^d [0,x_i]}=\chi_{[0,x]}$ for short:
\begin{eqnarray}
F &=&I_{1}(\sum_{i\geq 1}(h\chi_{[0,x]},e _{i})e _{i}(\xi)z)=\sum_{i\geq 1}(h\chi_{[0,x]},e
_{i})I_{1}(e _{i}(\xi)z) \\
&=&\sum_{i\geq 1}\left( \int_0^x h(\xi)e _{i}(\xi)d\xi\right)
K_{\varepsilon ^{(i,1)}}m_{2}.  \nonumber
\end{eqnarray}
Recall that $z=m_{2}p_{1}(z);$ see (\ref{10.5}).
In particular, choosing $h=1$ we get the following expansion for the L\' evy sheet $L$:
\begin{equation}
L(x)=\sum_{i\geq 1}\left( \int_{0}^{x}e _{i}(\xi)d\xi\right)
K_{\varepsilon ^{(i,1)}}m_{2}.  \label{LK}
\end{equation}

\subsection{The Skorohod integral}
The following definition of the {\it Skorohod integral} with respect to the L\'{e}vy field $L(x)$ is an extension to the  multiprameter case of a definition originally due to Kabanov
(1975):
\bigskip

\begin{definition}[The Skorohod integral]
\label{Skoro}
Let $Y(x)$, $x\in\mathbb{R}^d$, be a measurable stochastic process
such that
\[
\mathbb{E}[Y^2(x)]<\infty,
\qquad x\in\mathbb{R}^d.
\]
Then, for each $x\in\mathbb{R}^d$, $Y(x)$ admits a chaos expansion
of the form
\[
Y(x)=\sum_{n=0}^{\infty}I_n(f_n(\cdot,x)),
\]
where
\[
f_n(\cdot,x)\in
\widehat{L}^2\big((\lambda\times\nu)^n\big),
\qquad n\geq 1,
\]
and
\[
I_0(f_0(\cdot,x))=\mathbb{E}[Y(x)].
\]

For $n\geq 0$, let $\widetilde{f}_n$ denote the symmetrization, with
respect to all $n+1$ variables, of the function
\[
\begin{aligned}
&(x^{(1)},z_1,\ldots,x^{(n+1)},z_{n+1})
\\
&\qquad\longmapsto
z_{n+1}
f_n(x^{(1)},z_1,\ldots,x^{(n)},z_n,x^{(n+1)}).
\end{aligned}
\]

Suppose that
\[
\sum_{n=0}^{\infty}
(n+1)!
\|\widetilde{f}_n\|^2_{
L^2((\lambda\times\nu)^{n+1})}
<\infty.
\]
Then $Y(\cdot)$ is said to be \emph{Skorohod integrable} with respect
to $L$, and its Skorohod integral is defined by
\begin{equation}
\label{Skor}
\int_{\mathbb{R}^d}Y(x)L(dx)
:=
\sum_{n=0}^{\infty}
I_{n+1}(\widetilde{f}_n).
\end{equation}
\end{definition}

Note that 

$${\rm
\E}\left[\left(\int_{\R^d}Y(x)L(dx)\right)^2\right]=\sum_{n=0}^\infty(n+1)!||\tilde{f}_n||^2_{L^2((\lambda^d\times\nu)^{n+1})}<\infty$$
and therefore $\int_{\R^d}Y(x)L(dx)\in L^2(\P)$. Moreover, 

$$\E \left[\int_{\R^d}Y(x)L(dx)\right]=0.$$

\subsection{The L\'{e}vy-Hida spaces and the Kondratiev spaces}
We now introduce some useful stochastic test function spaces and stochastic distribution spaces.\\

In the following we put
\begin{equation}
(2\mathbb{N})^{k\alpha }=\prod_{j\geq 1}(2j)^{k\alpha _{j}},
\end{equation}%
if $\alpha =(\alpha _{1,}\alpha _{2,},...)\in \mathcal{J}$.

\begin{definition}
 (The L\'{e}vy-Hida-Kondratiev spaces of stochastic test functions and stochastic
distributions)

\parindent5mm
\item{{\bf a)}} {\it (The Kondratiev stochastic test function spaces)}. For $0\le\rho\le 1, k=1,2,...$, let $({\cal S})_{\rho,k}=({\cal
S})_{\rho,k}^{(L)}$ consist of those
$$\phi=\sum_{\alpha\in{\cal J}} c_\alpha K_\alpha(\omega)\in L^2(\P)\hskip1cm (c_\alpha\in\R\ \hbox{constants})$$
such that

$$||\phi||^2_{\rho,k}:=\sum_{\alpha\in{\cal J}}c_\alpha^2(\alpha!)^{1+\rho}(2\N)^{k\alpha}<\infty\hskip1cm \hbox{for {\it all} 
}k\in\N $$

\item{{\bf b)}} {\it (The Kondratiev stochastic distribution spaces)}. For $0\le\rho\le 1, k=1.2...$ let $({\cal S})_{-\rho,-k}=({\cal
S})_{-\rho,-k}^{(L)}$ consist of all formal expansions
$$F=\sum_{\alpha\in{\cal J}}b_\alpha K_\alpha(\omega)$$
such that

$$||F||^2_{-\rho,-k}:=\sum_{\alpha\in{\cal J}}b_\alpha^2(\alpha!)^{1-\rho}(2\N)^{-k\alpha}<\infty\hskip1cm \hbox{for {\it some}
}k\in\N$$

We define 
$$({\cal S})_\rho=\cap_{k=1}^{\infty} \S_{\rho,k},$$
equipped with the projective topology (intersection) defined by the norms $||\cdot||_{\rho,k};\
k=1,2,\dots$ and we define 
$$(\S)_{-\rho}=\cup_{k=1}^{\infty}(\S)_{-\rho,-k},$$ 
equipped with the inductive topology (union) defined by the norms 
$||\cdot||_{-\rho,-k};\ k=1,2,\dots$. Then $({\cal S})_{-\rho}$ becomes the dual of $({\cal S})_\rho$ and the action of $F=\sum
b_\alpha K_\alpha\in ({\cal S})_{-\rho}$ on $\phi=\sum a_\alpha K_\alpha\in ({\cal S})_{\rho}$ is given by

$$<F,\phi>=\sum_\alpha a_\alpha b_\alpha \alpha!\leqno{(5.3.30)}$$
\item{{\bf c)}} The space $({\cal S})=({\cal S})_0^{(L)}$ is called the {\it L\' evy-Hida stochastic test function space} and the space
$({\cal S})^*=({\cal S})_{-0}^{(L)}$ is called the {\it L\' evy-Hida distribution space}.\\

{\bf Generalized expectation.}\newline
If  $F= \sum_{\alpha \in \mathcal{J}} a_{\alpha }K_{\alpha } \in (\mathcal{S})^{* }$, we define
the generalized expectation $\E[F]$ of $F$ by
$$
\E[F] = a_0.
$$ 
Note that $\E[K_{\alpha }] = 0$ for all $\alpha \ne 0$.
Therefore the generalized expectation coincides with the usual expectation if $F \in L^2(\P)$.

\end{definition}

\subsection{White noise for L\' evy sheets and their compensated Poisson random measures}
Now we have the prerequisites to define the white noise of a L\'evy sheet and its associated Poisson random measure. In view of \eqref{LK} the following definition is natural:

\begin{definition} 
\begin{enumerate}
\item[(i)]
The \emph{L\'{e}vy white noise}
$$\overset{\bullet }{L
}(x):=V(x):= \frac{\partial^d L(x)}{\partial x_1 \partial x_2  ... ,\partial x_d},\quad x\in \R^d$$
is defined by the expansion
\begin{eqnarray}
V(x):&=&m_{2}\sum_{i\geq 1}e _{i}(x)K_{\varepsilon
^{(i,1)}}=m_{2}\sum_{i\geq 1}e _{i}(x)I_{1}(e _{i}(x)p_{1}(z))\nonumber\\
&=&m_{2}\sum_{i\geq 1}e _{i}(x)I_{1}(e _{i}(x)z).  
\end{eqnarray}%

\item[(ii)]
The \emph{white noise} $\overset{\bullet }{\widetilde{N}%
}(x,z)$ of the compensated random measure  $\widetilde{N}(dx,dz)$ is defined by the expansion%
\begin{equation}
\overset{\bullet }{\widetilde{N}}(x,z)=\sum_{i,j\geq 1}e
_{i}(x)p_{j}(z)K_{_{\varepsilon ^{(i,j)}}}(\omega ).  
\end{equation}
\end{enumerate}
\end{definition}

\begin{remark}
\label{Rem10.6}Note that in the case of the expansion for $\overset{\bullet }{L}(x); x \in \R^d$ we have%
\[
\sum_{\alpha \in \mathcal{J}}c_{\alpha }^{2}\alpha !(2\mathbb{N})^{-q\alpha
}=m_{2}^{2}\sum_{i\geq 1}e _{i}^{2}(x) 2^{-q}\kappa (i,1)^{-q}<\infty \text{ }
\]%

for $q\geq 2.$ Here we have used that $\kappa (i,1)=1+(i-1)i/2\geq i,$  and the following
well-known estimate for the Hermite functions:%
\begin{equation}
\sup_{x\in \mathbb{R}}\left\vert e _{m}(x)\right\vert =\mathcal{O}(m^{-\frac{1}{12}%
}).
\end{equation}%
This estimate is well-known for the 1-parameter Hermite functions with $x \in \R$, and it is easy to see that it extends to the tensor products $e_m(x)$ with $x \in \R^n$.
Therefore $\overset{\bullet }{L}(x)\in (\mathcal{S})^{* }$ for all $x$.
 Similarly $\overset{\bullet }{\widetilde{N}}(x,z)\in (\mathcal{S})^{* }$
for all $x,z.$
\end{remark}

\begin{remark}
\begin{enumerate}
\item[(i)]
By comparing the expansions%
\[
L(x)=m_2 \sum_{i\geq 1}\left( \int_{0}^{x_1}\int_0^{x_2} ... \int_0^{x_n} e _{i}(\xi) d\xi\right)
K_{\varepsilon ^{(i,1)}}
\]%
and%
\[
\overset{\bullet }{L}(x)=m_{2}\sum_{i\geq 1}e _{i}(x)K_{\varepsilon
^{(i,1)}},
\]%
we get formally
\begin{equation}
\overset{\bullet }{L}(x)=\frac{\partial^n L(x)}{\partial x_1 \partial x_2 ..., \partial x_n}\text{ (derivative in }(%
\mathcal{S})^{* }). 
\end{equation}
This can be proved rigorously.
\item[(ii)] 
Choose a Borel set $\mho $ such that $\overline{\mho }\subset \mathbb{%
R}\backslash \{0\}.$ Then%
\begin{eqnarray*}
\widetilde{N}(x,\mho ) &=&I_{1}\left( \chi _{[0,x]}(\xi)\chi _{\mho
}(z)\right)  \\
&=&\sum_{i,j\geq 1}(\chi _{[0,x]},e _{i})_{L^{2}(\mathbb{R}^n)}(\chi
_{\mho },p_{j})_{_{L^{2}(\nu )}}I_{1}\left( e _{i}(\xi)p_{j}(\zeta)\right) \\
&=&\sum_{i,j\geq 1}\int_{0}^{x}e _{i}(\xi)d\xi\cdot \int_{\mho }p_{j}(z)\nu
(z)K_{\varepsilon ^{(i,j)}}(\omega ).
\end{eqnarray*}%
This justifies the identity
\begin{equation}
\overset{\bullet }{\widetilde{N}}(x,z)=\frac{\widetilde{N}(dx,dz)}{dx\times
\nu (dz)}\text{ (Radon-Nikodym derivative).}  \label{10.26}
\end{equation}%
\item[(iii)]
Also note that $\overset{\bullet }{L}$ is related to $%
\overset{\bullet }{\widetilde{N}}$ by%
\begin{equation}
\overset{\bullet }{L}(x)=\int_{\realio} \overset{\bullet }{\widetilde{N%
}}(x,z)z\nu (dz).  
\end{equation}%
\dproof
To see this we consider%
\begin{eqnarray*}
&&\int_{\realio}\overset{\bullet }{\widetilde{N}}(x,z)z\nu (dz) \\
&=&\int_{\realio}\sum_{i,j\geq 1}e _{i}(x)p_{j}(z)K_{_{\varepsilon
^{(i,j)}}}(\omega )z\nu (dz)
=\sum_{i\geq 1}e _{i}(x)\sum_{j\geq 1}I_{1}( e
_{i}p_j)\int_{\realio} p_{j}(z)z\nu (dz))  \\
&=&m_2\sum_{i\geq 1}e _{i}(x)I_{1}\left( e _{i}(\xi)z\right) =\overset{%
\bullet }{L}(x).
\end{eqnarray*}
\fproof
\end{enumerate}
\end{remark}

\section{White noise calculus}
\subsection{The Wick product}

The Wick product is important in the study of
stochastic (ordinary and partial) differential equations. In
general, one can say that the use of this product corresponds to
-- and extends naturally -- the use of It\^{o} integrals. We now
explain this in more detail.

\begin{definition}The {\it Wick product} $F\diamond G$ of two elements
$$F=\sum\limits_\alpha a_\alpha K_\alpha,\;G=\sum\limits_{\beta} b_{\beta}
K_{\beta}\in({\cal S})^*\quad\text{with}\quad
a_{\alpha},b_{\beta}\in\R$$
is defined by
$$F\d 
G=\sum\limits_{\alpha,\beta}a_{\alpha}b_{\beta}K_{\alpha+\beta}.$$
\end{definition}
\begin{example}
 Choose 
 \begin{align}
 &F=L(x)=I_{1}(\sum_{i\geq 1}(\chi_{[0,x]},e _{i})e _{i}(\xi)z)\\
&=\sum_{i\geq 1}(\chi_{[0,x]},e
_{i})I_{1}(e _{i}(\xi)z)=\sum_{i\geq 1}\left( \int_0^x e _{i}(\xi)d\xi\right)
K_{\varepsilon ^{(i,1)}}m_{2}\\
&\text{ and }\\
&G=L(y)=I_{1}(\sum_{i\geq 1}(\chi_{[0,y]},e _{i})e _{i}(\xi)z)\\
&=\sum_{i\geq 1}(\chi_{[0,y]},e
_{i})I_{1}(e _{i}(\xi)z)=\sum_{i\geq 1}\left( \int_0^y e _{i}(\xi)d\xi\right)
K_{\varepsilon ^{(i,1)}}m_{2}.  
\end{align}
Then
\begin{align*}
    &L(x) \diamond L(y)= \sum_{i,j} ( \int_0^x e _{i}(\xi)d\xi)( \int_0^y e _{j}(\xi)d\xi) K_{\epsilon^{(i,1)}+\epsilon ^{(j,1)}}\\
    &=\sum_{i,j} ( \int_0^x e _{i}(\xi)d\xi)( \int_0^y e _{j}(\xi)d\xi) K_{\epsilon^{(i.1)+(j,1)}}\\
    &=\sum_{i,j} ( \int_0^x e _{i}(\xi)d\xi)( \int_0^y e _{j}(\xi)d\xi)I_2((e_ip_1)\widehat{\otimes} (e_j p_1)) \\
    &=\sum_{i,j} ( \int_0^x e _{i}(\xi)d\xi)( \int_0^y e _{j}(\xi)d\xi)\int_{\R}\int_{\R_0} \Big(\int_{\R} \int_{\R_0} e_i(x^{(1)})z_1e_j(x^{(2)})z_2 \widetilde{N}(dx^{(1)},dz_1)\Big) \widetilde{N}(dx^{(2)},dz_2)
\end{align*}

\end{example}

\subsection{Wick products and Skorohod integration}

We now prove a fundamental relation between Wick products and
Skorohod integrals:

\begin{definition}
\label{definition3.13-added} A function $Y:\mathbb{R}^d \to 
(\mathcal{S})^{*}$ (also called an $(\mathcal{S})^{*}$-valued process) is
\emph{$(\mathcal{S})^{*}$-integrable} if
\[
< Y(\cdot),f> \in L^{1}(\mathbb{R}^d), \quad\text{ for all}\quad f\in(\mathcal{S}).
\]
Then the \emph{$(\mathcal{S})^{*}$-integral of $Y$}%
\index{$(\mathcal{S})^{*}$-integral}%
, denoted by $\int_{\mathbb{R}^d} Y(x)dx$, is the (unique) element in
$(\mathcal{S})^{*}$ such that
\begin{equation}
\label{3.68-added}<\int_{\mathbb{R}^d}Y(x)dx, f> = \int_{\mathbb{R}^d} <Y(x),f>
dx, \quad\text{ for all } f\in(\mathcal{S}).
\end{equation}

\end{definition}

\begin{remark}
\label{remark3.14-added} The fact that \eqref{3.68-added} does indeed define
$\int_{\mathbb{R}^d }Y(x)dx$ as an element of $(\mathcal{S})^{*}$ is a
consequence of \cite[Proposition 8.1]{HKPS}.
\end{remark}

Note that if $Y$ is $(\mathcal{S})^{*}$-integrable, then so is $Y\chi_{(a,b]}%
$, for all $a,b\in\mathbb{R}^d$, and we put
\[
\int_{a}^{b} Y(x)dx := \int_{\mathbb{R}^d}Y(x)\chi_{(a,b]}(x)dx.
\]

The following result was first proved by Lindstr\o m et al \cite{LOU} (1992) and subsequently with a different proof by Benth \cite{B} (1993), both in the classical 1-parameter Brownian motion case. Our proof here is an adaptation of Benth's proof in \cite{B} to the multiparameter L\' evy sheet case.

\begin{theorem}[{\cite{LOU}, \cite{B}}]
\label{wick}
Let $Y:\mathbb{R}^d\times\Omega\to\mathbb{R}$ be a Skorohod integrable
stochastic process in the sense of Definition~\ref{Skoro}. Then the
$(\mathcal{S})^{*}$-valued process
\[
x\longmapsto Y(x)\diamond\overset{\bullet}{L}(x)
\]
is $(\mathcal{S})^{*}$-integrable and
\begin{equation}
\label{3.66}
\int_{\mathbb{R}^d}
Y(x)\diamond\overset{\bullet}{L}(x)\,dx
=
\int_{\mathbb{R}^d}Y(x)L(dx),
\end{equation}
where the integral on the right-hand side is the Skorohod integral
defined by
\[
\int_{\mathbb{R}^d}Y(x)L(dx)
:=
\sum_{n=0}^{\infty}
I_{n+1}\big(\widetilde{f}_n\big),
\]
with $\widetilde{f}_n$ denoting the symmetrization of
\[
z_{n+1}f_n
\]
in the $n+1$ variables.
\end{theorem}
\dproof
Let
\[
Y(x)=\sum_{\alpha\in\mathcal J}c_\alpha(x)K_\alpha
     =\sum_{n=0}^{\infty}I_n(f_n(\cdot,x)),
\]
where, for $|\alpha|=n$,
\begin{equation}
f_n(\cdot,x)
=
\sum_{|\alpha|=n}
c_\alpha(x)\,
\theta^{\hat\otimes\alpha}.
\label{fn-expansion}
\end{equation}
By assumption, $Y$ is Skorohod integrable, so that
\[
\sum_{n=0}^{\infty}
(n+1)!
\|\widetilde f_n\|^2_{L^2((\lambda\times\nu)^{n+1})}
<\infty,
\]
where $\widetilde f_n$ denotes the symmetrization of
\[
z_{n+1}f_n
(x^{(1)},z_1,\ldots,x^{(n)},z_n,x^{(n+1)}).
\]

We first compute the Wick integral. Recall that
\[
\overset{\bullet}L(x)
=
m_2\sum_{k=1}^{\infty}
e_k(x)K_{\varepsilon^{(k,1)}},
\]
where
\[
m_2=\left(\int_{\mathbb R_0}z^2\nu(dz)\right)^{1/2}.
\]
Using the definition of the Wick product,
\[
K_\alpha\diamond K_\beta=K_{\alpha+\beta},
\]
we obtain
\begin{align}
Y(x)\diamond\overset{\bullet} L(x)
&=
m_2
\sum_{\alpha\in\mathcal J}
\sum_{k=1}^{\infty}
c_\alpha(x)e_k(x)
K_{\alpha+\varepsilon^{(k,1)}}.
\end{align}
Consequently, by the definition of the $(\mathcal S)^*$-integral,
\begin{align}
\int_{\mathbb R^d}
Y(x)\diamond\overset{\bullet} L(x)\,dx
&=
m_2
\sum_{\alpha\in\mathcal J}
\sum_{k=1}^{\infty}
\left(
\int_{\mathbb R^d}
c_\alpha(x)e_k(x)\,dx
\right)
K_{\alpha+\varepsilon^{(k,1)}} \\
&=
m_2
\sum_{\alpha\in\mathcal J}
\sum_{k=1}^{\infty}
(c_\alpha,e_k)_{L^2(\mathbb R^d)}
K_{\alpha+\varepsilon^{(k,1)}}.
\label{wick-expansion}
\end{align}

We now compute the Skorohod integral using Definition~\ref{Skoro}.
For $|\alpha|=n$, we have
\[
f_n(\cdot,x)
=
\sum_{|\alpha|=n}
c_\alpha(x)\theta^{\hat\otimes\alpha}.
\]
Hence
\[
z_{n+1}f_n
=
\sum_{|\alpha|=n}
c_\alpha(x)
z_{n+1}
\theta^{\hat\otimes\alpha}.
\]
Since
\[
z=m_2p_1(z),
\]
and
\[
\theta_{\kappa(k,1)}(x,z)
=
e_k(x)p_1(z),
\]
we may expand
\[
c_\alpha(x)
=
\sum_{k=1}^{\infty}
(c_\alpha,e_k)_{L^2(\mathbb R^d)}e_k(x)
\]
and therefore obtain
\begin{align}
z_{n+1}f_n
&=
m_2
\sum_{|\alpha|=n}
\sum_{k=1}^{\infty}
(c_\alpha,e_k)_{L^2(\mathbb R^d)}
\,
\theta^{\hat\otimes\alpha}
\otimes
\theta_{\kappa(k,1)}.
\end{align}
After symmetrization,
\begin{align}
\widetilde f_n
&=
m_2
\sum_{|\alpha|=n}
\sum_{k=1}^{\infty}
(c_\alpha,e_k)_{L^2(\mathbb R^d)}
\,
\theta^{\hat\otimes\alpha}
\hat\otimes
\theta_{\kappa(k,1)}.
\end{align}
By the definition of the multi-index notation,
\[
\theta^{\hat\otimes\alpha}
\hat\otimes
\theta_{\kappa(k,1)}
=
\theta^{\hat\otimes
(\alpha+\varepsilon^{(k,1)})}.
\]
Thus
\begin{equation}
\widetilde f_n
=
m_2
\sum_{|\alpha|=n}
\sum_{k=1}^{\infty}
(c_\alpha,e_k)_{L^2(\mathbb R^d)}
\theta^{\hat\otimes
(\alpha+\varepsilon^{(k,1)})}.
\label{tilde-fn}
\end{equation}

By Definition~\ref{Skoro},
\begin{align}
\int_{\mathbb R^d}Y(x)L(dx)
&=
\sum_{n=0}^{\infty}
I_{n+1}(\widetilde f_n)\\
&=
m_2
\sum_{n=0}^{\infty}
\sum_{|\alpha|=n}
\sum_{k=1}^{\infty}
(c_\alpha,e_k)_{L^2(\mathbb R^d)}
I_{n+1}
\left(
\theta^{\hat\otimes
(\alpha+\varepsilon^{(k,1)})}
\right).
\end{align}
Since
\[
K_\beta
=
I_{|\beta|}
\left(
\theta^{\hat\otimes\beta}
\right),
\]
we have
\[
I_{n+1}
\left(
\theta^{\hat\otimes
(\alpha+\varepsilon^{(k,1)})}
\right)
=
K_{\alpha+\varepsilon^{(k,1)}}.
\]
Therefore
\begin{equation}
\int_{\mathbb R^d}Y(x)L(dx)
=
m_2
\sum_{\alpha\in\mathcal J}
\sum_{k=1}^{\infty}
(c_\alpha,e_k)_{L^2(\mathbb R^d)}
K_{\alpha+\varepsilon^{(k,1)}}.
\label{skor-expansion}
\end{equation}

Comparing \eqref{wick-expansion} and \eqref{skor-expansion}, we conclude that
\[
\int_{\mathbb R^d}
Y(x)\diamond\overset{\bullet} L(x)\,dx
=
\int_{\mathbb R^d}Y(x)L(dx).
\]
This proves the result.
\fproof
\begin{remark}
The identity in Theorem~\ref{wick} should be understood with
the Skorohod integral defined by the chaos expansion above. In the
one-parameter case, this integral coincides with the corresponding
It\^o integral for adapted integrands.
In the multiparameter setting, one has to distinguish carefully
between the chaos-defined Skorohod integral and the various
constructions of multiparameter It\^o integrals. The identity
proved above concerns the Skorohod integral defined by the chaos
expansion. The relation between this integral and the
multiparameter It\^o integral requires additional considerations.
The computation on page 229 of
\cite{HOUZ} concerns the latter integral and therefore does not
directly contradict Theorem~\ref{wick}.
\end{remark}
\begin{example}[The one-parameter case]
Let $d=1$ and consider the L\' evy process $L=(L(t))_{t\geq0}$.
We illustrate Theorem~\ref{wick} by taking
\[
Y(s)=L(s).
\]
By the Wick product calculus,
\[
\frac{d}{ds}L(s)^{\diamond 2}
=
2L(s)\diamond\overset{\bullet}{L}(s)
\]
in $(\mathcal S)^*$, and hence
\begin{equation}
\label{wick-square}
\int_0^t
L(s)\diamond\overset{\bullet}{L}(s)\,ds
=
\frac12L(t)^{\diamond2}.
\end{equation}

On the other hand, Theorem~\ref{wick} gives
\[
\int_0^t
L(s)\diamond\overset{\bullet}{L}(s)\,ds
=
\int_0^t L(s)L(ds),
\]
where the integral on the right-hand side is the Skorohod integral
defined through the chaos expansion.

Since $L(s)$ is adapted, in the one-parameter case the Skorohod
integral coincides with the It\^o integral. Thus
\[
\int_0^t L(s)L(ds)
=
\int_0^t L(s-)L(ds).
\]
Applying the It\^o formula for the square of a pure-jump L\' evy process,
we obtain
\begin{equation}
\label{ito-square}
L(t)^2
=
2\int_0^t L(s^-)L(ds)
+
\int_0^t\int_{\mathbb R_0}z^2N(ds,dz).
\end{equation}
Combining \eqref{wick-square} and \eqref{ito-square}, we obtain
\[
L(t)^{\diamond2}
=
L(t)^2
-
\int_0^t\int_{\mathbb R_0}z^2N(ds,dz).
\]
Equivalently,
\[
L(t)^{\diamond2}
=
L(t)^2-[L]_t,
\]
where $[L]_t$ denotes the quadratic variation of $L$.
\end{example}
\begin{example}
We illustrate Theorem \ref{wick} in the two-parameter case.
Let
\[
L(x)=L(x_1,x_2), \qquad x=(x_1,x_2)\in\mathbb{R}_+^2,
\]
be the L\' evy sheet and set
\[
Y(x)=L(x).
\]

Recall that the L\' evy sheet admits the representation
\[
L(x)
=
\int_{[0,x]}\int_{\mathbb{R}_0}
z\,\widetilde{N}(dy,dz),
\]
where
\[
[0,x]=[0,x_1]\times[0,x_2].
\]

Define the deterministic kernel
\[
\phi_x(y,z)
:=
\mathbf{1}_{[0,x]}(y)\,z,
\qquad
(y,z)\in\mathbb{R}^2\times\mathbb{R}_0.
\]
Consequently,
\[
Y(x)=L(x)=I_1(\phi_x).
\]
Hence, the process \(Y\) is Skorohod integrable
whenever the corresponding integrability condition in Definition
\ref{Skoro} is satisfied. Hence, by Theorem \ref{wick},
\[
\int_{[0,x]}L(y)\diamond\overset{\bullet} L(y)\,dy
=
\int_{[0,x]}L(y)\,L(dy),
\]
where the integral on the right-hand side is the Skorohod integral
defined by Definition \ref{Skoro}.

This identity should be distinguished from the multiparameter It\^o
integral. In particular, for \(d=2\), the Wick-calculus identity
\[
\frac12 L(x)^{\diamond2}
=
\int_{[0,x]}L(y)\diamond\overset{\bullet} L(y)\,dy
+
\int_{[0,x]}
\frac{\partial L}{\partial y_1}(y)
\diamond
\frac{\partial L}{\partial y_2}(y)\,dy
\]
shows that, in general,
\[
\frac12 L(x)^{\diamond2}
\neq
\int_{[0,x]}L(y)\,L(dy).
\]
There is therefore no contradiction with Theorem \ref{wick}: the
additional term is a genuinely multiparameter phenomenon and reflects
the distinction between the Skorohod integral and the multiparameter
It\^o integral. Moreover, to the best of our knowledge, a general
It\^o formula for multiparameter L\'evy sheets is not available.
\end{example}

\vspace{2mm} In view of the above theorem, the following terminology is natural.

\begin{definition}
\label{Wick-definition6.16} Suppose $Y$ is an $(\mathcal{S})^{\ast}$-valued
process such that
\[
\int_{\mathbb{R}^d}Y(x)\diamond\overset{\bullet}{L}(x)dx:\in(\mathcal{S})^{\ast
},
\]
then we call this integral \emph{the generalized Skorohod integral of $Y(x)$} with respect to $L(dx)$.
\end{definition}

Combining the properties above we get a powerful calculation technique for Skorohod 
integration. First of all, note that, by Theorem \ref{wick},
\begin{equation}
\int\limits_{0}^{x}\overset{\bullet}{L}(\xi)d\xi=L({x})\,. \label{3.63}%
\end{equation}
One can deduce that
\begin{equation}
\int\limits_{0}^{x}Z\diamond Y({\xi})\diamond\overset{\bullet}{L}({\xi})d\xi
=Z\diamond\int\limits_{0}^{x}Y({\xi})\diamond\overset{\bullet}{L}({\xi})d\xi,
\label{3.64}%
\end{equation}
if $Z$ does not depend on $\xi$. This is in contrast to the fact that for Skorohod
integrals we generally have
\begin{equation}
\int\limits_{0}^{x}Z\cdot Y({\xi})L(d\xi))\not =Z\cdot\int\limits_{0}%
^{x}Y({\xi})L(d\xi)\,, \label{3.65}%
\end{equation}
even if $Z$ does not depend on $\xi$.

\section{The population growth equation}
We illustrate the theory developed above by applying it to study the following model for population growth in a random environment driven by both multiparameter Brownian and  L\' evy white noise:

\begin{align}\label{heatgrowth}
 	 \frac{\partial \partial^{\alpha}}{\partial x\partial t^{\alpha}}Y(t,x)&= \mu(t,x) Y(t,x)+\sigma(t,x)Y(t,x) \diamond W(t,x)+\gamma(t,x)Y(t,x)\diamond V(t,x);\nonumber\\
     & (t,x)\in (0,\infty)\times \mathbb{R}^{d}.
 	 	\end{align}
with initial conditions

\begin{equation}
Y(0,x) = y_0, \quad Y(t,0) = y_0 \quad \text{(constant) for all } t,x \ge 0,
\end{equation}
where:
\begin{itemize}
    \item $\frac{\partial^\alpha}{\partial t^\alpha}$ is the Caputo fractional derivative of order $\alpha \in (0,1)$,
    \item $W(t,x) = \frac{\partial^2 B(t,x)}{\partial t \partial x}$ is the Brownian white noise (Brownian sheet),
    \item $V(t,x) = \frac{\partial^2 L(t,x)}{\partial t \partial x}$ is the L\'evy white noise (pure jump L\'evy sheet),
    \item $\mu(t,x)$, $\sigma(t,x)$, $\gamma(t,x)$ are deterministic functions,
    \item $y_0$ is a constant initial population density.
\end{itemize}
\begin{remark}
    Note that if $\alpha=1$ this corresponds to the following SDE driven by Brownian and L\' evy sheets:
\begin{align}\label{heatgrowth}
 	 Y(t,x)&= \int_0^t \int_0^x \mu(s,a) Y(s,a)ds da+\int_0^t \int_0^x \sigma(s,a)Y(s,a) B(ds,da)\nonumber\\
     &+\int_0^t \int_0^x \gamma(s,a)Y(s,a)L(ds,da);\nonumber\\
     & (t,x)\in (0,\infty)\times \mathbb{R}^{d}.
 	 	\end{align}
        \end{remark}
        
Before we state and prove our main result, we recall the following useful estimate on the norm of Wick products in the Kondratiev spaces:

\begin{proposition}(The V\r{a}ge inequality \cite{HOUZ}, Proposition 3.3.2)\\
Suppose $F\in (\S)_{-1,-k}$ and $G \in (\S)_{-1,-\ell}$ for some integers $k>\ell +1$.
Then
\begin{align}
 \|F \diamond G \|_{-1,-k} \leq A(k-\ell) \|F\|_{-1,-k} \|G\|_{-1,-\ell},  
\end{align}
where
\begin{align}
    A(k-\ell) = \sum_{\alpha} (2\mathbb{N})^{(\ell -k)\alpha}.
\end{align}
\end{proposition}
From this we deduce the following by induction on $n$:
\begin{corollary} \label{5.2}
    Suppose $K_i \in (\S)_{-1,-\ell}$ for some $\ell > 2; i=1,2, ...$
    Then for $k> \ell +1$ we have
    \begin{align}
        \|K_1 \diamond K_2 \diamond \cdots \diamond K_n\|_{-1,-k} \leq A(k-\ell)^{n-1} \Pi_{i=1}^n \|K_i\|_{-1,-\ell}
    \end{align}
    In particular, if $\|K_i\|_{-1,-\ell} \leq C$ for all $i$ then
  \begin{align}
        \|K_1 \diamond K_2 \diamond \cdots \diamond K_n\|_{-1,-k} \leq A(k-l)^{n-1} C^n.
    \end{align}  
\end{corollary}

We now state and prove our main result:
\begin{theorem}
    The solution $Y(t,x)$ of the fractional stochastic population growth equation \ref{heatgrowth} exists and is unique as an element of the Hida stochastic distribution space $(\S)^{*}$ and is explicitly represented as
\begin{align}
    Y(t,x)&=y_0 \sum_{n=0}^\infty \frac{1}{n!n!} \Big(\int_{[0,u]} K(v) dv \Big)^{\d n}
\end{align}
where \\
$K(v)=K(s,z)=\frac{1}{\Gamma(\alpha)}(t-s)^{\alpha -1}\Big(\mu(s,z) +\sigma(s,z)\overset{\bullet}{B}(s,z)  + \gamma(s,z) \overset{\bullet}{L}(s,z)\Big)$\\
and hence
\begin{align}
 &\int_{[0,u]} K(v) dv\nonumber\\
 &=\int_0^t \int_0^x \frac{1}{\Gamma(\alpha)}(t-s)^{\alpha -1}\Big(\mu(s,z) +\sigma(s,z)\overset{\bullet}{B}(s,z)  + \gamma(s,z) \overset{\bullet}{L}(s,z)\Big)ds dz\nonumber \\
 &= \int_0^t\int_0^x \frac{1}{\Gamma(\alpha)}(t-s)^{\alpha -1}\mu(s,z)ds dz\nonumber \\
 &+\int_0^t\int_0^x  \frac{1}{\Gamma(\alpha)}(t-s)^{\alpha -1}\sigma(s,z) B(ds,dz) +\int_0^t\int_0^x  \frac{1}{\Gamma(\alpha)}(t-s)^{\alpha -1}\gamma(s,z) L(ds,dz).
\end{align}    
\end{theorem}

\dproof
(i) We first prove existence and uniqueness:\\
By integrating  with respect to the spatial variable from $0$ to $x$:

\begin{align}
\frac{\partial^\alpha Y(t,x)}{\partial t^\alpha} - \frac{\partial^\alpha Y(t,0)}{\partial t^\alpha} &= \int_0^x \mu(t,z) Y(t,z) \, dz \notag \\
&\quad + \int_0^x \sigma(t,z) Y(t,z) \diamond W(t,z) \, dz \notag \\
&\quad + \int_0^x \gamma(t,z) Y(t,z) \diamond V(t,z) \, dz.
\end{align}

Since $Y(t,0) = y_0$, we obtain:

\begin{equation}\label{fracderiv}
\frac{\partial^\alpha Y(t,x)}{\partial t^\alpha} = \int_0^x \mu(t,z) Y(t,z) \, dz + \int_0^x \sigma(t,z) Y(t,z) \diamond W(t,z) \, dz + \int_0^x \gamma(t,z) Y(t,z) \diamond V(t,z) \, dz.
\end{equation}

Now apply the Caputo fractional integral of order $\alpha$ to both sides. Recall that the Caputo fractional integral $I^\alpha$ satisfies:
 \begin{equation}
I^\alpha \left( \frac{\partial^\alpha Y}{\partial t^\alpha} \right)(t,x) = Y(t,x) - Y(0,x).
\end{equation}
 
Hence:

\begin{align}
Y(t,x) &= y_0 + \frac{1}{\Gamma(\alpha)} \int_0^t (t-s)^{\alpha-1} \int_0^x \mu(s,z) Y(s,z) \, dz \, ds \notag \\
&\quad + \frac{1}{\Gamma(\alpha)} \int_0^t (t-s)^{\alpha-1} \int_0^x \sigma(s,z) Y(s,z) \diamond W(s,z) \, dz \, ds \notag \\
&\quad + \frac{1}{\Gamma(\alpha)} \int_0^t (t-s)^{\alpha-1} \int_0^x \gamma(s,z) Y(s,z) \diamond V(s,z) \, dz \, ds.
\end{align}   
That can be written as follow:

\begin{align}\label{integraleq}
Y(t,x) &= y_0 + \frac{1}{\Gamma(\alpha)} \int_0^t (t-s)^{\alpha-1} \int_0^x \mu(s,z) Y(s,z) \, dz \, ds \notag \\
&\quad + \frac{1}{\Gamma(\alpha)} \int_0^t (t-s)^{\alpha-1} \int_0^x \sigma(s,z) Y(s,z) \diamond B(ds, dz) \notag \\
&\quad + \frac{1}{\Gamma(\alpha)} \int_0^t (t-s)^{\alpha-1} \int_0^x \gamma(s,z) Y(s,z) \diamond L(ds, dz).
\end{align}

To solve this equation, we use the $S$-Transform method:

Let $F$ be a stochastic distribution in the Hida space. The S-transform of $F$ is defined by

\[
S(F)(\varphi) = \left\langle F , :e^{\langle \omega, \varphi \rangle}:\right\rangle.
\]

It satisfies

\[
S(F \diamond G)(\varphi) = S(F)(\varphi) S(G)(\varphi).
\]

We apply the S transform  to the equation \eqref{integraleq}:

Define
\[
u(t,x,\varphi) := S(Y(t,x))(\varphi).
\]

We have:

\[
S(\mu Y) = \mu(s,z)\ u(s,z,\varphi),
\]

\[
S(Y \diamond B(ds,dz)) = u(s,z,\varphi)\  \varphi(s,z)\ ds\ dz,
\]

\[
S(Y \diamond L(ds,dz)) = u(s,z,\varphi)\ \psi(s,z)\ ds\ dz,
\]

where
\[
\psi(s,z) = \int_{\mathbb{R}} (e^{\varphi(s,z,u)} - 1)\  \nu(du).
\]

Using all these, we obtain the following deterministic equation:

\begin{equation}
\begin{aligned}
u(t,x,\varphi)
&= y_0 \\
&+ \frac{1}{\Gamma(\alpha)} \int_0^t (t-s)^{\alpha-1}
\int_0^x \mu(s,z)\ u(s,z,\varphi)\  dz\  ds \\
&+ \frac{1}{\Gamma(\alpha)} \int_0^t (t-s)^{\alpha-1}
\int_0^x \sigma(s,z)\ u(s,z,\varphi)\ \varphi(s,z)\  dz\  ds \\
&+ \frac{1}{\Gamma(\alpha)} \int_0^t (t-s)^{\alpha-1}
\int_0^x \gamma(s,z)\ u(s,z,\varphi)\ \psi(s,z)\  dz\  ds.
\end{aligned}
\end{equation}

Define
\[
\beta_\varphi(s,z) = \mu(s,z) + \sigma(s,z)\ \varphi(s,z) + \gamma(s,z)\ \psi(s,z).
\]

Then
\begin{equation}\label{volt}
  u(t,x,\varphi) = y_0 + \frac{1}{\Gamma(\alpha)} \int_0^t (t-s)^{\alpha-1} \int_0^x \beta_\varphi(s,z)\ u(s,z,\varphi) dz  ds.  
\end{equation}

The S-transform reduces the stochastic equation to a deterministic linear fractional Volterra-type integral equation.

For a fixed test function $\varphi$, define the Volterra operator

\[
(K_\varphi f)(t,x)
=
\frac{1}{\Gamma(\alpha)}
\int_0^t (t-s)^{\alpha-1}
\int_0^x
\beta_\varphi(s,z)f(s,z)\,dz\,ds.
\]
Then equation \eqref{volt} can be written as

\[
u=y_0+K_\varphi u.
\]

We define the Picard sequence by
\[
u_0(t,x,\varphi)=y_0,
\]

and

\[
u_{n+1}(t,x,\varphi)
=
y_0+K_\varphi u_n(t,x,\varphi).
\]
The first iterates are

\[
u_1
=
y_0+K_\varphi y_0,
\]

\[
u_2
=
y_0+K_\varphi y_0+K_\varphi^2y_0,
\]

\[
u_3
=
y_0+K_\varphi y_0+K_\varphi^2y_0+K_\varphi^3y_0.
\]
By induction,

\[
u_n(t,x,\varphi)
=
\sum_{k=0}^{n}
K_\varphi^k y_0.
\]
We now compute $K_\varphi^k y_0$.\\ For $k=1$,

\[
K_\varphi y_0
=
\frac{y_0}{\Gamma(\alpha)}
\int_0^t
(t-s_1)^{\alpha-1}
\int_0^x
\beta_\varphi(s_1,z_1)\,dz_1\,ds_1.
\]
 
 For $k=2$,

\[
K_\varphi^2 y_0
=
\frac{y_0}{\Gamma(\alpha)^2}
\int_{0<s_2<s_1<t}
(t-s_1)^{\alpha-1}
(s_1-s_2)^{\alpha-1}
\]
\[
\qquad\qquad\times
\int_{0<z_2<z_1<x}
\beta_\varphi(s_1,z_1)
\beta_\varphi(s_2,z_2)
\,dz_2\,dz_1\,ds_2\,ds_1.
\]
Proceeding recursively, one obtains

\[
K_\varphi^k y_0
=
\frac{y_0}{\Gamma(\alpha)^k}
\int_{0<s_k<\cdots<s_1<t}
\prod_{j=1}^{k}
(s_{j-1}-s_j)^{\alpha-1}
\]

\[
\qquad\qquad\times
\int_{0<z_k<\cdots<z_1<x}
\prod_{j=1}^{k}
\beta_\varphi(s_j,z_j)
\,dz_k\cdots dz_1\,
ds_k\cdots ds_1,
\]

where $s_0=t$.
Assume that for every fixed $\varphi$ and every compact $[0,T]\times[0,X]$,
\[
B_\varphi(T,X)
:=
\sup_{(t,x)\in[0,T]\times[0,X]}
|\beta_\varphi(t,x)| < \infty.
\]
Then for any bounded function $f$,
\begin{align}
|K_\varphi f(t,x)|
&\le
\frac{B_\varphi(T,X)}{\Gamma(\alpha)}
\int_0^t (t-s)^{\alpha-1}
\int_0^x |f(s,z,\varphi)|\,dz\,ds.
\end{align}
Using the definition of the iterated integral in $\mathbb{R}^m$,
\[
\int_0^x dz = \prod_{i=1}^m x_i,
\quad \text{when } f \equiv 1,
\]
we obtain the bound
\begin{align}
\int_0^x |f(s,z)|\,dz
\le
\Big(\prod_{i=1}^m x_i\Big)\|f\|_\infty.
\end{align}
Let
\[
X=(X_1,\ldots,X_m)\in\mathbb{R}_+^m,
\qquad
[0,X]:=\prod_{i=1}^m[0,X_i].
\]
For \(x=(x_1,\ldots,x_m)\in[0,X]\), we have
\[
0\leq x_i\leq X_i,\qquad i=1,\ldots,m.
\]
Hence\begin{align}
\|K_\varphi f\|_\infty
\le
\frac{B_\varphi(T,X)}{\Gamma(\alpha)}
\Big(\prod_{i=1}^m X_i\Big)
\int_0^t (t-s)^{\alpha-1} ds\;\|f\|_\infty.
\end{align}
 Since
\[
\int_0^t (t-s)^{\alpha-1} ds
=
\frac{t^\alpha}{\alpha}
=
\frac{t^\alpha}{\Gamma(\alpha+1)}\Gamma(\alpha),
\]
we get
\begin{align}
\|K_\varphi f\|_\infty
\le
C_{\varphi,T,X}\, t^\alpha \|f\|_\infty,
\end{align}
where
\[
C_{\varphi,T,X}
=
\frac{B_\varphi(T,X)\prod_{i=1}^m X_i}{\Gamma(\alpha+1)}.
\]
By induction,
\[
\|K_\varphi^k y_0\|_\infty
\le
|y_0|
\frac{\big(C_{\varphi,T,X} T^\alpha\big)^k}{\Gamma(k\alpha+1)}.
\]
Since $\Gamma(k\alpha+1)$ grows super-exponentially, the series
\[
u(t,x,\varphi)
=
\sum_{k=0}^{\infty} K_\varphi^k y_0
\]
converges uniformly on $[0,T]\times[0,X]$.
Using linearity of $K_\varphi$,
\[
K_\varphi u
=
\sum_{k=1}^{\infty} K_\varphi^k y_0,
\]
hence
\[
y_0 + K_\varphi u
=
\sum_{k=0}^{\infty} K_\varphi^k y_0
=
u.
\]
For each fixed $\varphi \in \mathcal{S}(\mathbb{R}^d)$, there exists a unique solution
\[
u(t,x,\varphi)
\in C([0,T]\times[0,X])
\]
given by the Picard-Neumann series
\begin{equation}
    u(t,x,\varphi)
=
\sum_{k=0}^{\infty} K_\varphi^k y_0.
\end{equation}

\
It remains to reconstruct the stochastic solution from its
$S$-transform. For each fixed $(t,x)$, consider the mapping
\[
\varphi \longmapsto u(t,x,\varphi).
\]

This mapping is a
\emph{$U$-functional}. More precisely, it verifies the following
\begin{enumerate}
\item for all $\varphi_1,\varphi_2 \in \mathcal{S}(\mathbb{R}^d)$, the function
\[
z \longmapsto u(t,x,z\varphi_1+\varphi_2),
\qquad z\in\mathbb{C},
\]
is entire;

\item there exist constants $C,a,p>0$ such that
\[
|u(t,x,z\varphi)|
\le
C\exp\!\big(a|z|^2\|\varphi\|_p^2\big),
\]
for all $z\in\mathbb{C}$ and $\varphi\in\mathcal{S}(\mathbb{R}^d)$.
\end{enumerate}

The characterization
theorem for the $S$-transform in the Hida-Kondratiev framework
 implies that there exists a unique
stochastic distribution
\[
Y(t,x)\in (\mathcal{S})^{*}
\]
such that
\[
S(Y(t,x))(\varphi)
=
u(t,x,\varphi),
\qquad
\varphi\in\mathcal{S}(\mathbb{R}^d).
\]

Therefore, $u(t,x,\varphi)$ is the $S$-transform of a unique stochastic
distribution $Y(t,x)$, and we define
\[
Y(t,x)
=
S^{-1}\big(u(t,x,\cdot)\big).
\]

Consequently, $Y(t,x)$ is the unique stochastic distribution solution
of the original stochastic fractional integral equation.\\

(ii)
Next, to get an explicit solution formula, we work directly in the Kondratiev spaces $(\S)_{-1,-r}$ using Wick calculus: \\

Note that the equation can be written
\begin{align}
Y(t,x)=y_0 + \int_0^t \int_0^x K(s,z) \diamond Y(s,z) ds dz, \label{4.1}
\end{align}
where
\begin{align}
K(s,z)=K_t(s,z)=\frac{1}{\Gamma(\alpha)}(t-s)^{\alpha -1}(\mu(s,z) +\sigma(s,z)\overset{\bullet}{B}(s,z)  + \gamma(s,z) \overset{\bullet}{L}(s,z)).\label{4.2}
\end{align}
Substituting for $Y(s,z)$ in \eqref{4.1} we get
\begin{align*}
Y(u)= y_0 + \int_0^{u} K(u_1) \d \big(y_0+ \int_0^{u_1}K(u_2)\d Y(u_2)du_2\big)du_1,
\end{align*}
where we put  $u_j=(s_j,z_j); j=1,2, ...$.\\ 
Repeating this we obtain
\begin{align}
&Y(u)=y_0\Big( 1 + \int_0^{u} K(u_1)du_1 + \int_0^{u} (\int_0^{u_1} K(u_1) \d K(u_2)du_2 ) du_1 \nonumber\\
& +  ... + \int_0^{u}  \int_0^{u_1} ... \int_0^{u_n} K(u_1)\d K(u_2) \d \cdots K(u_{n}) du_1 du_2 \cdots du_n) +R_{n} \label{4.5}\\
&=y_0
+y_0
\sum_{k=1}^{n}
\int_{\Delta_k(u)}
K(u_1)\diamond K(u_2)\diamond\cdots\diamond K(u_k)
\,du_k\cdots du_1
+
R_n(u)
\end{align}
where
\[
\Delta_k(u)
=
\{0<s_k<\cdots<s_1<t,\; 0<z_k<\cdots<z_1<x\}.
\]
and 
\begin{align}
R_{n}&= \int_0^{u} \big( \int_0^{u_1}\cdots \int_0^{u_n} K(u_1)\d K(u_2) \d ... K(u_{n+1}) \d Y(u_{n+1})du_1 du_2 ...du_n \big) du_{n+1}\\
&=\int_{\Delta_{n+1}(u)}
K(u_1)\diamond\cdots\diamond K(u_{n+1})
\diamond Y(u_{n+1})
\,du_{n+1}\cdots du_1.
\end{align}
By the V\r{a}ge inequality (Corollary \ref{5.2}), it follows that $R_{n} \rightarrow 0$ in $(\mathcal{S})_{-1,-r}$ as $n \rightarrow \infty$.

Define 
$$F(u_1,u_2,\cdots, u_n)=K(u_1)\diamond K(u_2) \diamond \cdots. \diamond K(u_n).$$
Then the function $F$ is symmetric in all its variables $u_1, u_2,\cdots, u_n$.
This is simply by the fact that $\diamond$ is a commutative multiplication.
\\

Therefore we get, by letting $n \rightarrow \infty$ in \eqref{4.5},
\begin{align}
&Y(t,x)= y_0 \sum_{n=0}^\infty \int_0^{u} \big( \int_0^{u_1} ... \int_0^{u_n} K(u_1)\d K(u_2) \d ... K(u_n) du_1 du_2 ...du_n\\
&=y_0
\sum_{n=0}^{\infty}
\int_{\Delta_n(u)}
K(u_1)\diamond K(u_2)\diamond\cdots\diamond K(u_{n})
\,du_{1}\cdots du_n\\
&=y_0 \sum_{n=0}^\infty \frac{1}{n! n!} \int_{[0,u]^n} K(u_1)\d K(u_2) \d ... \d K(u_n) du_1 ... du_n\\
&=y_0 \sum_{n=0}^\infty \frac{1}{n!n!} \Big(\int_{[0,u]} K(v) dv \Big)^{\d n} \label{Y}
\end{align}
where
\begin{align}
 &\int_{[0,u]} K(v) dv\nonumber\\
 &=\int_0^t \int_0^x \frac{1}{\Gamma(\alpha)}(t-s)^{\alpha -1}\Big(\mu(s,z) +\sigma(s,z)\overset{\bullet}{B}(s,z)  + \gamma(s,z) \overset{\bullet}{L}(s,z)\Big)ds dz\nonumber \\
 &= \int_0^t\int_0^x \frac{1}{\Gamma(\alpha)}(t-s)^{\alpha -1}\mu(s,z)ds dz\nonumber \\
 &+\int_0^t\int_0^x  \frac{1}{\Gamma(\alpha)}(t-s)^{\alpha -1}\sigma(s,z) B(ds,dz) +\int_0^t\int_0^x  \frac{1}{\Gamma(\alpha)}(t-s)^{\alpha -1}\gamma(s,z) L(ds,dz).
\end{align}
The first part of the proof shows that the solution is unique and belongs to  $(\S)^{*}$. Therefore we can conclude that the solution $Y(t,x)$ given by \eqref{Y} actually belongs to $(\S)^{*}$.
\fproof

\begin{remark}
    Observe that the solution obtained in \eqref{Y} is not the Wick exponential. Indeed, the Wick exponential is defined by
\[
\exp^{\diamond}(Z)
=
\sum_{n=0}^{\infty}
\frac{1}{n!}
Z^{\diamond n},
\]
whereas our solution is given by
\[
Y(t,x)
=
y_0
\sum_{n=0}^{\infty}
\frac{1}{(n!)^2}
Z(t,x)^{\diamond n},
\]
with
\[
Z(t,x)
=
\int_0^u K(v)\,dv.
\]
Recall that the modified Bessel function of the first kind of order zero is defined by the power series
\[
I_0(z)
=
\sum_{n=0}^{\infty}
\frac{1}{(n!)^2}
\left(\frac{z}{2}\right)^{2n},
\]
or equivalently,
\[
I_0(2\sqrt{z})
=
\sum_{n=0}^{\infty}
\frac{z^n}{(n!)^2}.
\]
Therefore, the series appearing in \eqref{Y} is precisely the power series of the modified Bessel function evaluated at the Wick variable \(Z\). This motivates the introduction of the Wick--Bessel function
\[
I_0^{\diamond}(2\sqrt{Z})
:=
\sum_{n=0}^{\infty}
\frac{Z^{\diamond n}}{(n!)^2},
\]
where the ordinary powers are replaced by Wick powers.
Consequently, the solution can be written formally as
\begin{equation}
    Y(t,x)
=
y_0
I_0^{\diamond}
\!\left(
2\sqrt{
\int_0^u K(v)\,dv
}
\right).
\end{equation}
\end{remark}

\section{Acknowledgments}
 This research was funded by the Research Council of Norway under project 363480 (NASTRAN-Talent).

\end{document}